\documentclass{gtart_h}

\def\ifplaintex{\expandafter\ifx\csname documentclass\endcsname\relax}

\def\gtp{{\mathsurround=0pt\it $\cal G\mskip-2mu$eometry \&\ 
$\cal T\!\!$opology $\cal P\!$ublications}}  

\def\recd{{\small Received:\qua\receiveddate\ifx\reviseddate\relax
\else\qquad Revised:\qua\reviseddate\fi\par}} 

\def\lognumber#1{\def\thelognumber{#1}}
\def\volumenumber#1{\def\thevolumenumber{#1}}
\def\volumeyear#1{\def\thevolumeyear{#1}}
\def\papernumber#1{\def\thepapernumber{#1}}
\def\pagenumbers#1#2{\def\startpage{#1}\def\finishpage{#2}}
\def\published#1{\def\publishdate{#1}}

\def\received#1{\def\receiveddate{#1}}
\def\revised#1{\def\reviseddate{#1}}
\def\accepted#1{\def\accepteddate{#1}}

\def\asciiurl#1{\def\theasciiurl{#1}}

\let\\\par\let\thelognumber\relax\let\thevolumenumber\relax
\let\thepapernumber\relax\let\thevolumeyear\relax\let\startpage\relax
\let\finishpage\relax\let\publishdate\relax\let\receiveddate\relax
\let\reviseddate\relax\let\accepteddate\relax\let\theasciititle\relax
\let\theasciiauthors\relax
\let\theasciiabstract\relax

\let\theasciiemail\relax
\let\theasciiurl\relax

\ifplaintex
\font\logobig=cmssbx10 scaled 3836
\font\logomed=cmssbx10 scaled 2557
\else
\font\logobig=cmssbx10 scaled 4200
\font\logomed=cmssbx10 scaled 2800
\fi

\long\def\makeagttitle{   
\count0=\startpage
\agt\hfill      
\hbox to 45truept{\vbox to 0pt{\vglue -13truept{\logomed A\kern -.37em{\logobig 
T}\kern -.38em G}\vss}\hss}
\break
{\small Volume \thevolumenumber\ (\thevolumeyear)
\startpage--\finishpage\nl
Published: \publishdate}

\vglue .25truein

{\parskip=0pt\leftskip 0pt plus
1fil\def\\{\par\smallskip}{\Large\bf\thetitle}\par\medskip} \vglue
0.05truein

{\parskip=0pt\leftskip 0pt plus 1fil\def\\{\par}{\sc\theauthors}
\par\medskip}%
 
\vglue 0.03truein 

{\small\leftskip 25truept\rightskip 25truept{\bf Abstract}\stdspace\theabstract

{\bf AMS Classification}\stdspace\theprimaryclass
\ifx\thesecondaryclass\relax\else; \thesecondaryclass\fi\par
{\bf Keywords}\stdspace \thekeywords\par}\vglue 7truept

}   

\ifplaintex
\font\phead=cmsl9 scaled 950
\font\pnum=cmbx10 scaled 913
\font\pfoot=cmsl9 scaled 950
\headline{\vbox to 0pt{\vskip -4.5mm\line{\small\phead\ifnum
\count0=\startpage ISSN 1472-2739 (on-line) 1472-2747 (printed)
\hfill {\pnum\folio}\else\ifodd\count0\def\\{ }%
\ifx\theshorttitle\relax\thetitle\else\theshorttitle\fi\hfill{\pnum\folio}
\else\def\\{ and }{\pnum\folio}\hfill\ifx\theshortauthors\relax\theauthors
\else\theshortauthors\fi\fi\fi}\vss}}
\footline{\vbox to 0pt{\vglue 0mm\line{\small\pfoot\ifnum\count0=\startpage
\copyright\ \gtp\hfill\else
\agt, Volume \thevolumenumber\ (\thevolumeyear)\hfill\fi}\vss}}
\else
\font=cmsl9 scaled 1050
\font\lnum=cmbx10 
\font=cmsl9 scaled 1050
\makeatletter
\def\@oddhead{{\small\lhead\ifnum\count0=\startpage ISSN 1472-2739 
(on-line) 1472-2747 (printed)\hfill {\lnum\number\count0}\else\ifodd\count0
\def\\{ }\ifx\theshorttitle\relax \thetitle \else\theshorttitle\fi\hfill
{\lnum\number\count0}\else\def\\{ and }{\lnum\number\count0}
\hfill\ifx\theshortauthors\relax 
\theauthors\else\theshortauthors\fi\fi\fi}}\def\@evenhead{\@oddhead}
\def\@oddfoot{\small\lfoot\ifnum\count0=\startpage\copyright\ \gtp\hfill\else
\agt, Volume \thevolumenumber\ (\thevolumeyear)\hfill\fi}
\def\@evenfoot{\@oddfoot}
\makeatother
\fi
\let\maketitlepage\makeagttitle
\let\makeshorttitle\maketitlepage
\let\maketitle\maketitlepage

\newwrite\gtoutfile
\long\gdef\makeheadfile{  
{\def\\{, }\def\s{ }
\immediate\openout\gtoutfile head.xxx
\immediate\write\gtoutfile{Proxy-for: \ifx\theasciiauthors\relax
\theauthors\else\theasciiauthors\fi\s<\ifx\theasciiemail\relax\theemail\else\theasciiemail\fi>}
\immediate\write\gtoutfile{\noexpand\\}
\immediate\write\gtoutfile{Authors: \ifx\theasciiauthors\relax
\theauthors\else\theasciiauthors\fi}
{\def\\{ }\immediate\write\gtoutfile{Title: \ifx\theasciititle\relax
\thetitle\else\theasciititle\fi}}
\immediate\write\gtoutfile{Subj-class: GT or SG, GR etc}
\immediate\write\gtoutfile{MSC-class: \theprimaryclass\ifx\thesecondaryclass\relax\else, \thesecondaryclass\fi}
\immediate\write\gtoutfile{Journal-ref: Algebraic and Geometric Topology \thevolumenumber\s
(\thevolumeyear) \startpage-\finishpage}
\immediate\write\gtoutfile{Comments: Published by Algebraic and
Geometric Topology at}
\immediate\write\gtoutfile{\s\s\s  http://www.maths.warwick.ac.uk/agt/AGTVol\thevolumenumber/agt-\thevolumenumber-\thepapernumber.abs.html}
\immediate\write\gtoutfile{\noexpand\\}
\immediate\write\gtoutfile{}
\ifx\theasciiabstract\relax
\immediate\write\gtoutfile{\theabstract}\else
\immediate\write\gtoutfile{\theasciiabstract}\fi
\immediate\write\gtoutfile{}
\immediate\write\gtoutfile{\noexpand\\}
\immediate\write\gtoutfile{}
\immediate\closeout\gtoutfile}}  

\def\maketitlepage{\makeagttitle\makeheadfile}
\let\makeshorttitle\maketitlepage
\let\maketitle\maketitlepage

\lognumber{13}
\volumenumber{4}
\volumeyear{2004}
\papernumber{13}
\published{13 April 2004}
\pagenumbers{219}{241}
\received{23 October 2003}
\revised{20 January 2004}
\accepted{23 January 2004}

\usepackage{amsmath, amssymb}
\newtheorem{theorem}{Theorem}[section]

\newtheorem{corollary}[theorem]{Corollary}

\newtheorem{lemma}[theorem]{Lemma}
\newtheorem{proposition}[theorem]{Proposition}
\theoremstyle{definition}

\newtheorem{definition}[theorem]{Definition}
\newtheorem*{remark}{Remark}

\begin{document}

\title{Adem relations in the Dyer-Lashof algebra\\and modular invariants}
\shorttitle{Adem relations in the Dyer-Lashof algebra}
\authors{Nondas E. Kechagias}
\address{Department of Mathematics, University of Ioannina, 45110 Greece}
\email{nkechag@cc.uoi.gr}
\gturl{\url{http://www.math.uoi.gr/~nondas_k}}
\asciiurl{http://www.math.uoi.gr/ nondas k}

\begin{abstract}
This work deals with Adem relations in the Dyer-Lashof algebra from a
modular invariant point of view.  The main result is to provide an
algorithm which has two effects: Firstly, to calculate the hom-dual of
an element in the Dyer-Lashof algebra; and secondly, to find the image
of a non-admissible element after applying Adem relations. The
advantage of our method is that one has to deal with polynomials
instead of homology operations. A moderate explanation of the
complexity of Adem relations is given.
\end{abstract}

\maketitle

\section{Introduction}

The relationship between the (canonical sub-co-algebras) Dyer-Lashof
algebra, $R[k]$ and the Dickson invariants $D[k]$ is well-known, see May's
paper in \cite{C-L-M}, relevant parts of which will be quoted here. We
provide an algorithm for calculating Adem relations in the Dyer-Lashof
algebra using modular co-invariants. Much of our work involves the
calculation of the hom-duals of elements of $R$ in terms of the generators
of the polynomial algebra $D[k]$. The results described here will be applied
to give an invariant theoretic description of the mod$-p$ cohomology
of a finite loop space in \cite{kech4}.

We note that the idea for our algorithm was inspired by May's theorem 3.7,
page 29, in \cite{C-L-M}. The key ingredient for relating homology
operations and polynomial invariants is the relation between the map which
imposes Adem relations and the decomposition map between certain rings of
invariants. This relation was studied by Mui for $p=2$ in \cite{Mui2}, and
we extend it here for any prime. Namely:

\medskip
\textbf{Theorem 4.15 }\qua\textsl{Let }$\rho :T[n]\rightarrow R[n]$\textsl{\ be
the map which imposes Adem relations. Let }$\hat{\imath}:S(E(n))^{GL_{n}}%
\otimes D[n]\hookrightarrow S(E(n))^{B_{n}}\otimes B[n]$\textsl{\ be the
natural inclusion. Then }$\rho ^{\ast }\equiv \hat{\imath}$\textsl{, i.e.\
for any }$e_{I}\in T[n]$\textsl{\ and }$d^{m}M^{\varepsilon }\in
S(E(n))^{GL_{n}}\otimes D[n]$\textsl{, } 
$$\langle d^{m}M^{\varepsilon },\rho (e_{I})\rangle =\langle \hat{\imath}(d^{m}M^{\varepsilon }),e_{I}\rangle .
$$

Campbell, Peterson and Selick studied self maps $f$ of $\Omega
_{0}^{m+1}S^{m+1}$ and proved that if $f$ induces an isomorphism on $%
H_{2p-3}(\Omega _{0}^{m+1}S^{m+1},\mathbb Z /p \mathbb Z)$, then $f_{(p)}$ is a homotopy
equivalence for $p$ odd and $m$ even \cite{C-P-S}. A key ingredient for
their proof was the calculation of 
$$AnnPH^{\ast }(\Omega _{0}^{m+1}S^{m+1},\mathbb Z /p \mathbb Z)$$
They gave a convenient method for calculating the hom-dual of elements of $%
H_{\ast }(\Omega _{0}^{m+1}S^{m+1},\mathbb Z /p \mathbb Z)$ which do not involve Bockstein
operations. Our algorithm computes the hom-duals of elements of $R[n]$ in
terms of the generators of the polynomial algebra $D[n]$. Please see Theorem 
\ref{dual}.

A direct application of the last two theorems is the computation of Adem
relations. The main difference between the classical and our approach is
that we consider Adem relations ``globally'' instead of consecutive elements
and it requires fewer calculations. This algorithm is described in
Proposition \ref{Adem-Rel}.

The paper is purely algebraic and its applications are deferred to \cite
{kech4}. There are three sections in this paper beyond this introduction,
sections 2, 3 and 4. Section 2 recalls well known facts about the
Dyer-Lashof algebra from May's article, cited above. In section 3, the
Dickson algebra and its relation with the ring of invariants of the Borel
subgroup is examined. That relation is studied using a certain family of
matrices which suitably summarizes the expressions for Dickson invariants in
terms of the invariants of the Borel subgroup. In the view of the author,
the complexity of Adem relations is reflected in the different ways in which
the same monomial in the generators of the Borel subgroup can show up as a
term in a Dickson invariant. The ways in which this can happen can be
understood using these matrices. For $p$ odd, the dual of the Dyer-Lashof
algebra is a subalgebra of the full ring of invariants. This subalgebra is
also discussed in full details. In the last section a great amount of work
is devoted to the proof of the analog of Mui's result mentioned above. Then
our algorithms more or less naturally follows.

This paper has been written for odd primes with minor modifications needed
when $p=2$ provided in statements in square brackets following the odd
primary statements.

For the sake of accessibility we shorten proofs. A detailed version of this
work including many examples can be found at: \url{http://www.uoi.gr/~nondas_k}
and also at:\newline
\url{http://www.maths.warwick.ac.uk/agt/ftp/aux/agt-4-13/full.ps.gz}

This work is dedicated to the memory of Professor F.P. Peterson.

We thank Eddy Campbell very much for his great effort regarding the
presentation and organization of the present work and the referee for his
encourangment and valuable suggestions regarding the accessibility of our
algorithm to the interested reader. Last but not least, we thank the editor
very much.

\section{The Dyer-Lashof algebra}

Let us briefly recall the construction of the Dyer-Lashof algebra. Let $F$
be the free graded associative algebra on $\{f^{i},\ i\geq0\}$ and $\{\beta
f^{i},\ i>0\}$ over $K:=\mathbb Z /p \mathbb Z$ with $|f^{i}|=2(p-1)i$, [$|f^{i}|=i$] and $%
|\beta f^{i}|=2i(p-1)-1$. $F$ becomes a co-algebra equipped with coproduct $%
\psi:\ F\longrightarrow F\otimes F$ given by 
$$\psi f^{i}=\sum f^{i-j}\otimes f^{j}\ \text{and}\ \psi\beta f^{i}=\sum\beta
f^{i-j}\otimes f^{j}\ +\ \sum f^{i-j}\otimes\beta f^{j}.$$
Elements of $F$ are of the form $$f^{I,\varepsilon }=\beta ^{\epsilon
_{1}}f^{i_{1}}\dots \beta ^{\epsilon _{n}}f^{i_{n}}$$ where $(I,\varepsilon
)=((i_{1},\dots ,i_{n}),$ $(\epsilon _{1},...,\epsilon _{n}))$ with $\epsilon
_{j}=0\ $or$\ 1$ and $i_{j}$ a non-negative integer for $j=1,\ \dots \ ,n$, $%
|f^{I,\varepsilon }|=2(p-1)\left( \sum\limits_{t=1}^{n}i_{t}\right) -\left(
\sum\limits_{t=1}^{n}e_{t}\right) $ [$|f^{I,\varepsilon }|=\left(
\sum\limits_{t=1}^{n}i_{t}\right) $]. Let $l(I,\varepsilon )=n$ denote the
length of $I,\varepsilon $ or $f^{I,\varepsilon }$ and let the excess of $%
(I,\varepsilon )$ or $f^{I,\varepsilon }$ be denoted $exc(f^{I,\varepsilon
})=i_{1}-\epsilon _{1}-|f^{I_{2}}|$, where $(I_{t},\varepsilon
_{t})=((i_{t},\dots ,i_{n}),(\epsilon _{t},\dots ,\epsilon _{n}))$. 
$$exc(f^{I,\varepsilon })=i_{1}-\epsilon _{1}-2(p-1)\sum\limits_{2}^{n}i_{t}%
\text{, [}exc(f^{I})=i_{1}-\sum\limits_{2}^{n}i_{t}\text{]}$$
The excess is defined $\infty $, if $I=\emptyset $ and we omit the sequence $%
(\epsilon _{1},...,\epsilon _{n})$, if all $e_{i}=0$. We refer to elements $%
f^{I}$ as having non-negative excess, if $exc(f^{I_{t}})$ is non-negative
for all $t$.

It is sometimes convenient to use lower notation for elements of $F$ and its
quotients. We define $f^{i}x=f_{\frac{1}{2}(2i-|x|)}x$ [$f^{i}x=f_{i-|x|}x$%
]. Let $I=(i_{1},...,i_{n})$ and $\varepsilon =(\epsilon _{1},...,\epsilon
_{n})$, then the degree of $Q_{I,\varepsilon }$ is 
$$|f_{I,\varepsilon }|=2(p-1)\left( \sum\limits_{t=1}^{n}i_{t}p^{t-1}\right)
-\left( \sum\limits_{t=1}^{n}e_{t}p^{t-1}\right) \text{, [}%
|f_{I,\varepsilon }|=\left( \sum\limits_{t=1}^{n}i_{t}2^{t-1}\right) \text{]%
}.$$
In lower notation we see immediately that $f_{I,\varepsilon }$ has
non-negative excess if and only if $(I,\varepsilon )$ is a sequence of
non-negative integers: $exc(I,\varepsilon )=2i_{1}-e_{1}$.

Given sequences $I$ and $I^{\prime}$ we call the direct sum of $I$ and $%
I^{\prime}$ the sequence $I\oplus I^{\prime}=(i_{1},...,i_{n},i_{1}^{\prime
},...,i_{m}^{\prime})$. Using a sequence $I $ we use the above idea for the
appropriate decomposition. Let $0_{k}$ denote the zero sequence of length $k$%
.

\begin{remark}
Let $\langle \mathbb N ,\frac{1}{2}\rangle $ be the monoid generated by $\mathbb N$ and $\frac{1}{2}$ in
the rationals. Let $\langle \mathbb N ,\frac{1}{2}\rangle ^{n}$ be the monoid which is the $n$-th
Cartesian product of $\langle \mathbb N ,\frac{1}{2}\rangle $. Then $(I,\varepsilon )\in \langle \mathbb N ,\frac{1%
}{2}\rangle ^{n}\times (Z/2Z)^{n}$. [$I\in \mathbb N ^{n}$]
\end{remark}

$F$ admits a Hopf algebra structure with unit $\eta:K\longrightarrow F$ and
augmentation $\epsilon:F\longrightarrow K$ given by: 
\begin{equation*}
\epsilon(f^{i})=\left\{ 
\begin{tabular}{ll}
$1,$ & if $i=0$ \\ 
$0,$ & otherwise.
\end{tabular}
\right.
\end{equation*}

\begin{definition}
There is a natural order on the elements  $f_{(I,\varepsilon )}$ defined as
follows: for $(I,\varepsilon )$ and $(I^{\prime },\varepsilon ^{\prime })$
we say that $(I,\varepsilon )<(I^{\prime },\varepsilon ^{\prime })$ if $%
exc(I_{l},\varepsilon _{l})=exc(I_{l}^{\prime },\varepsilon _{l}^{\prime })$
for $1\leq l\leq t$ and $exc(I_{t},\varepsilon _{t})<exc(I_{t}^{\prime
},\varepsilon _{t}^{\prime })$ for some $1\leq t\leq n$.
\end{definition}

We define $T=F/\mathcal{I}_{exc}$, where $\mathcal{I}_{exc}$ is the two
sided ideal generated by elements of negative excess. $T$ inherits the
structure of a Hopf algebra and if we let $T[n]$ denote the set of all
elements of $T$ with length $n$, then $T[n]$ is a co-algebra of finite type.
We denote the image of $f_{I,\varepsilon }$ by $e_{I,\varepsilon }$. Degree,
excess and ordering for upper or lower notation described above passes to $T$
and $T[n]$.

The Adem relations are given by: 
\begin{equation}
e_{r}e_{s}=\sum\limits_{i}(-1)^{\mathbf{r}-i}\binom{(p-1)(i-s)-1}{r-i-1}%
e_{r+ps-pi}e_{i}\text{, if }r>s  \notag
\end{equation}
and if $p>2$ and $r\geq s$, 
\begin{align*}
e_{r}\beta e_{s}& =\sum\limits_{i}(-1)^{r+i+1/2}\binom{(p-1)(i-s)}{r-1/2-i}%
\beta e_{r+ps-pi-1/2}e_{i}+ \\
& \sum\limits_{i}(-1)^{r+i-1/2}\binom{(p-1)(i-s)-1}{r-1/2-i}e_{r+ps-pi}\beta
e_{i}.
\end{align*}

Let $\mathcal{I}_{Adem}$ be the two sided ideal of $T$ generated by the Adem
relations. We denote $R$ the quotient $T/\mathcal{I}_{Adem}$ and this
quotient algebra is called \textit{the Dyer-Lashof algebra.} $R$\textbf{\ }%
is a Hopf algebra and $R[n]$ is again a co-algebra of finite type. We will
denote the obvious epimorphism above which imposes Adem relations by 
\begin{equation*}
\rho :T\rightarrow R\text{ with }\rho (e_{I})=\sum a_{I,J}Q_{J}
\end{equation*}
If $(I,\varepsilon )$ is admissible then $Q_{I,\varepsilon }$ is the image
of $e_{I,\varepsilon }$. 

The following lemma will be applied in section 4.

\begin{lemma}
\label{111}{\rm a)}\qua $\rho(e_{p^{k}}e_{0})=Q_{0}Q_{p^{k-1}}$; $\rho(e_{1}e_{0})=0$.

{\rm b)}\qua $\rho(e_{p^{k}+1/2}e_{1/2})=Q_{1/2}Q_{p^{k-1}+1/2}$; $%
\rho(e_{3/2}e_{1/2})=0$.

{\rm c)}\qua $\rho(e_{p^{k}+1}e_{1})=Q_{1}Q_{p^{k-1}+1}$; $\rho(e_{2}e_{1})=0$.

{\rm d)}\qua $\rho(e_{p^{k}}e_{1})=Q_{0}Q_{p^{k-1}+1}$; $\rho(e_{p}e_{1})=2Q_{0}Q_{2} $%
.

{\rm e)}\qua $\rho(e_{p^{k}}\beta e_{1/2})=Q_{0}\beta Q_{p^{k-1}+1/2}$; $%
\rho(e_{1}\beta e_{1/2})=\beta Q_{1/2}Q_{1/2}$.

{\rm f)}\qua $\rho(e_{p^{k}+1/2}e_{1/2})=0$; $\rho(e_{3/2}\beta e_{1})=\beta
Q_{1}Q_{1} $.

{\rm g)}\qua $\rho(e_{p^{k}+1}\beta e_{1/2})=\beta Q_{1/2}Q_{p^{k-1}+1/2}$; $\rho
(e_{2}\beta e_{1/2})=0$.
\end{lemma}

The passage from lower to upper notation between elements of $R$ is given as
follows. Let  $Jx\varepsilon $ and $Ix\varepsilon $ be lower and upper
sequences as defined above. Then, 
$$\beta ^{\epsilon _{1}}Q_{j_{1}}...\beta ^{\epsilon _{n}}Q_{j_{n}}\equiv
\beta ^{\epsilon _{1}}Q^{i_{1}}...\beta ^{\epsilon _{n}}Q^{i_{n}}$$
up to a unit in $\mathbb Z /p \mathbb Z$, where $i_{n}=j_{n}$, and 
$$i_{n-t} =\frac{1}{2}(2j_{n-t}+|I_{n-t+1}x\varepsilon _{n-t+1}|),  
j_{n-t} =\frac{1}{2}(2i_{n-t}-|J_{n-t+1}x\varepsilon _{n-t+1}|)$$
\begin{definition}
We say that  an element $Q_{I,\varepsilon }$ is admissible, if $0\leq
2i_{t}-2i_{t-1}+e_{t-1}$ for $2\leq t\leq n-1$.
\end{definition}

The ordering described above passes to $R$ and $R[n]$.

Since $R[n]$ and $T[n]$ are of finite type, they are isomorphic to their
duals as vector spaces and these duals become algebras. We shall describe
these duals giving an invariant theoretic description, namely: they are
isomorphic to subalgebras of rings of invariants over the appropriate
subgroup of $GL(n,K)$ in section 4.

\section{The Dickson algebra and a special family of matrices}

The Dickson algebra is a universal object in modular invariants of finite
groups. Applications involve computations of Dickson invariants of different
height. We provide formulas of this nature which will be applied in the
proof of Theorem \ref{dual}. Being very technical, those formulas can be
studied easier using matrices. 

Let $V^{k}$ denote a $K$-dimensional vector space generated by $%
\{e_{1},...,e_{k}\}$ for $1\leq k\leq n$. Let the dual basis of $V^{n}$ be $%
\{x_{1},...,x_{n}\}$ and the contragradient representation of $W_{\Sigma
_{p^{n}}}(V^{n})\longrightarrow Aut(V^{n})\equiv GL_{n}$ induces an action
of $GL_{n}$ on the graded algebra $E(x_{1},...,x_{n})\otimes
P[y_{1},...,y_{n}]$, [$P[y_{1},...,y_{n}]$], where $\beta x_{i}=y_{i}$. Let $%
E(n)=E(x_{1},...,x_{n})$ and $S[n]=K\left[ y_{1},\cdots ,y_{n}\right] $. The
degree is given by $|x_{i}|=1$ and $|y_{i}|=2$ (if $p=2$, then $|y_{i}|=1$).

The following theorems are well known:

\begin{theorem}
{\rm\cite{Dic}}\qua $S[n]^{GL_{n}}:=D[n]=K\left[ d_{n,0},\cdots ,d_{n,n-1}\right] $,
the Dickson algebra, is a polynomial algebra and their degrees are $\left|
d_{n,i}\right| =2\left( p^{n}-p^{i}\right) $, [$2^{n}-2^{i}$].
\end{theorem}

\begin{theorem}
{\rm\cite{Mui}}\qua $S[n]:=B[n]=K\left[ h_{1},\cdots,h_{n}\right] $ is a polynomial
algebra and their degrees are $\left| h_{i}\right| =2p^{i-1}\left(
p-1\right) $, [$2^{i-1}$].
\end{theorem}

Although relations between generators of the last two algebras can be easily
described, it is not the case between invariants of parabolic subgroups of
the general linear group.

Let $f_{k-1}\left( x\right) =\prod\limits_{u\in V^{k-1}}\left( x-u\right) $,
then $f_{k-1}\left( x\right) =\sum\limits_{i=0}^{k-1}\left( -1\right)
^{n-i}x^{p^{i}}d_{k-1,i}$ and $h_{k}=\smash{\prod\limits_{u\in V^{k-1}}}\left(
y_{k}-u\right) $. Moreover, (see \cite{Kech}), 
\begin{equation}
d_{n,n-i}=\sum\limits_{1\leq j_{1}<\cdots<j_{i}\leq
n}\smash{\prod\limits_{s=1}^{i}}\left( h_{j_{s}}\right) ^{p^{n-i+s-j_{s}}}
\label{formula}
\end{equation}
Let $m=(m_{0},...,m_{n-1})$ and $k=(k_{1},...,k_{n})$ be sequences of
non-negative integers. Let $d^{m}$ denote an element of $D[n]$ given by $%
\prod\limits_{t=0}^{n-1}d_{n,t}^{m_{t}}$\ and $h^{k}$ denote an element of $%
B[n]$ given by $\prod\limits_{t=1}^{n}h_{t}^{k_{t}}$. Let $I_{(t)}$ denote
the $t$-th element of the sequence $I=(i_{l_{1}},...,i_{l_{n}})$ from the
left: i.e.\ $I_{(t)}:=i_{l_{t}}$.

For any non-negative matrix $C$ with integral entries and $\mathbf{1}%
=(1,...,1)$, the matrix product $\mathbf{1\cdot}C$ is a sequence of
non-negative integers, then $h^{\mathbf{1\cdot}C}$ stands for $\prod
\limits_{t=1}^{n}h_{t}^{(\mathbf{1\cdot}C\mathbf{)}_{(t)}}$ Let $%
C(d_{n,j})=\{h^{I}\in B[n]$ and $h^{I}$ is a non-trivial summand in $%
d_{n,j}\}$, then $C(d_{n,i})\cap C(d_{n,j})=\emptyset$ for $j\neq i$.

\begin{remark}
1)\qua \label{zero-matrix}Before we start considering sets of matrices, we would
like to stress the point that the zero matrix is excluded from our sets,
unless otherwise stated.

2)\qua Until the end of this section, we number matrices beginning with $(0,0)$
in the upper left corner. In this case $h^{\mathbf{1\cdot}C}$ stands for $%
\prod\limits_{t=1}^{n}h_{t}^{(\mathbf{1\cdot}C\mathbf{)}_{(t-1)}}$.
\end{remark}

Let $0\leq j\leq n-1$. Here $j$ corresponds to the value $n-i$ in formula 
\ref{formula}.

\begin{definition}
For each matrix $A=(a_{it})$ such that $a_{it}$ is a non-negative integer, $\sum\limits_{t=0}^{n-1}a_{jt}=n-j$ and $\sum\limits_{t=0}^{n-1}a_{it}=0$
for $i\neq j$, we define an $n\times n$ matrix $C(A)=(b_{ij})=\left(
b_{(0)},\cdots ,b_{(n-1)}\right) $ such that $b_{it}=a_{it}p^{i-1-t+a_{i0}+%
\cdots +a_{it}}$\label{lem1}. Let us call this collection $A_{n,j}$.
\end{definition}

For $C\in A_{n,j}$, $\mathbf{1\cdot}C$ is the $j$-th row of $C$ which is the
only non-zero row of that matrix.

Let us also note that there is an obvious bijection between $A_{n,j}$ and $%
C(d_{n,j})$.

\begin{lemma}
$d_{n,j}=\sum\limits_{C\in A_{n,j}}h^{\mathbf{1\cdot}C}$.
\end{lemma}

\begin{definition}
Let $m=(m_{0},\cdots ,m_{n-1})$ be a sequence of zeros or powers of $p$. Let 
$A_{n,j}^{m}=\{m\cdot C_{j}=(m_{0}b_{(0)},\cdots
,m_{n-1}b_{(n-1)})\;|\;C_{j}=(b_{(0)},\cdots ,b_{(n-1)})\in A_{n,j}\}$ and  $%
A_{n}^{m}=\{\sum\limits_{j=0}^{n-1}m\cdot C_{j}\;|\;C_{j}\in A_{n,j}\}$. 
\end{definition}

Note that different elements of $A_{n}^{m}$ may provide the same element of $%
B[n]$ and this is the reason why  Adem relations are complicated as we shall
examine more in Proposition \ref{Prop9}. We shall also note that the
motivation of this section was exactly to demonstrate this difficulty using
an elementary method.

The following lemma is easily deduced from formulae \ref{formula}.

\begin{lemma}
\label{+++}Let $m=(m_{0},\cdots,m_{n-1})$ such that $m_{i}=0$ or $p^{k_{i}}$%
, then 
\begin{equation*}
d^{m}=\prod\limits_{0}^{n-1}d_{n,i}^{m_{i}}=\sum\limits_{C\in
A_{n}^{m}}\prod\limits_{t=1}^{n}\left( h_{t}\right) ^{(C(\mathbf{1}))_{t-1}}.
\end{equation*}
\end{lemma}

Coefficients might appear in the last summation. Hence one needs to
partition the set $A_{n}^{m}$ as the following lemma suggests.

\begin{lemma}
Let $m=(m_{0},\cdots ,m_{n-1})$ be a sequence of zeros or powers of $p$. Let 
$A=(a_{it})$ and $A^{\prime }=(a_{it}^{\prime })$ such that $%
a_{it},a_{it}^{\prime }\in \mathbb N$, $\sum\limits_{t=0}^{n-1}a_{jt}=\sum%
\limits_{t=0}^{n-1}a_{jt}^{\prime }=n-j$ if $m_{i}\neq 0$, otherwise the
last sums are zero. Suppose that $\mathbf{1}\cdot A=\mathbf{1}\cdot
A^{\prime }$ and let $\{i_{1},...,i_{q}\}$ denote their different columns.
Consider only their different rows and for each column $i_{r}$ partition
them according to where $1$'s appear: $\{j_{1},...,j_{s}\}$ and $%
\{j_{1}^{\prime },...,j_{s}^{\prime }\}$. If for each $j_{t}$ there exists a 
$j_{\ell }^{\prime }$ such that the number of zeros next to $a_{i_{r},j_{t}}$
and $a_{i_{r},j_{t}^{\prime }}^{\prime }$ are equal and this is true for all 
$i_{r}$, then $\mathbf{1\cdot }C(A)=\mathbf{1\cdot }C(A^{\prime })$.
\end{lemma}

\proof%
%
We use the definition of $C(A)$ in \ref{lem1}.%
\endproof%
%

On $1xn$ or $nx1$ matrices we give the left or upper lexicographical
ordering respectively.

\begin{definition}
Let $m$ be a non-negative integer, we denote by $_{|A_{n,j}|}(m)$ the set
of partitions of $m$ in $|A_{n,j}|$ terms. A typical element of $%
_{|A_{n,j}|}(m)$ is of the form $\pi=(\pi_{1},...,\pi_{|A_{n,j}|})$.

For $\pi=(\pi_{1},...,\pi_{|A_{n,j}|})\in _{|A_{n,j}|}(m)$, let $(\pi)$
denote the integer $\frac{m!}{\prod \pi_{t}!}$.
\end{definition}

\begin{lemma}
Let $m_{j}=\sum \limits_{\alpha=0}^{\ell_{j}}m_{j,\alpha}p^{\alpha}$. Then 
$$d_{nj}^{m_{j}}=\sum\limits_{\substack{ 0\leq\alpha\leq\ell_{j} \\ \pi^{(j,\alpha)}\in _{|A_{n,j}|}(m_{j,\alpha})}}\prod\limits_{%
\alpha=0}^{\ell_{j}}\left( \pi^{(j,\alpha)}\right) h^{\sum
\limits_{\alpha}\pi_{i}^{(j,\alpha)}\sum\limits_{C_{j,i}\in
A_{n,j}}p^{\alpha}\mathbf{1\cdot }C_{j,i}}$$
\end{lemma}

\proof%
%
First, we show the formulae above for $m_{j,\alpha }$ and then we extend by
direct multiplication.%
\endproof%
%

\begin{proposition}
\label{Prop9}Let $m=(m_{0},...,m_{n-1})$ be a sequence of non-negative
integers, then 
$$d^{m}=\sum\limits_{\substack{ 0\leq j\leq n-1,0\leq\alpha\leq\ell_{j} \\ \pi ^{(j,\alpha)}\in _{|A_{n,j}|}(m_{j,\alpha})}}\prod\limits_{j=0}^{n-1}\prod \limits_{\alpha=0}^{\ell_{j}}\left(
\pi^{(j,\alpha)}\right) h^{\sum \limits_{j}\sum\limits_{\alpha}\pi_{i}^{(j,\alpha)}\sum\limits_{\substack{C_{j,i}\in A_{n,j} \\ 0\leq j\leq n-1}}p^{\alpha}\mathbf{1\cdot }C_{j,i}}$$
\end{proposition}

The following lemma which is of great importance for dealing with Adem
relations involving Bockstein operations is proved using appropriate
matrices.

\begin{lemma}
\label{18}Each term of $d_{k+t,s}$ is also a term of $d_{k,s}d_{k+t,k}$.
Here $0\leq s<k$ and $1\leq t$. Moreover, no term of $%
d_{k,s}d_{k+t,k}-d_{k+t,s}$ is divisible by $\prod \limits_{k+1}^{k+t}h_{i}$%
.
\end{lemma}

In order to prove the main theorem in the next section, the following
formula for decomposing Dickson generators will be needed. This formula is a
special case of the lemma above. Formulas of this kind might be of interest
for other circumstances involving the Dickson algebra. One of them may be
the transfer between the Dickson algebra and the ring of invariants of
parabolic subgroups.

\begin{lemma}
Let $0\leq s<k$. Then $d_{k,s}d_{k+1,k}-d_{k+1,s}=$\newline
$\sum
\limits_{t=0}^{s-1}d_{k-t-1,s-t}^{p^{t}}d_{k-t-1,k-t-2}^{p^{t+2}}h_{k-t}^{p^{t}}+d_{k-s,0}^{p^{s}}d_{k-s,k-s-1}^{p^{s+1}}+d_{k-s,1}^{p^{s-1}}h_{k-s+1}^{p^{s-1}} 
$.
\end{lemma}

\begin{proof}
We shall use induction and the well known formula $%
d_{k,s}=d_{k-1,s-1}^{p}+d_{k-1,s}h_{k}$. 
\end{proof}

\begin{lemma}
\label{20}Each term of $d_{k+q,k}d_{k+t,s}$ is also a term of $%
d_{k+q,s}d_{k+t,k}$. Here $0\leq s<k$ and $0\leq q<t$. Moreover, no term of $%
d_{k+q,s}d_{k+t,k}-d_{k+q,k}d_{k+t,s}$ is divisible by $\prod
\limits_{k+q+1}^{k+t}h_{i}$.
\end{lemma}

\begin{proof}
We consider $(k+t)\times(k+t)$ matrices of the following form: 
\begin{equation*}
\begin{array}{c}
\\ 
s\text{-th} \\ 
\\ 
k\text{-th}
\end{array}
\left[ 
\begin{array}{cccc}
&  & \overset{k+q}{|} &  \\ 
\leftarrow & { k+q-s\rightarrow} & | &  \\ 
&  & | &  \\ 
\leftarrow & t & | & \rightarrow \\ 
&  & | & 
\end{array}
\right] - 
\begin{array}{c}
\\ 
s\text{-th} \\ 
\\ 
k\text{-th}
\end{array}
\left[ 
\begin{array}{cccc}
&  & \overset{k+q}{|} &  \\ 
\leftarrow & { k+t-s} & | & \rightarrow \\ 
&  & | &  \\ 
\leftarrow & q\rightarrow & | &  \\ 
&  & | & 
\end{array}
\right]
\end{equation*}
The last column of the matrices above is of size $t-q$. If this column is
full of non-zero elements in the last matrix, we require the same in the $k$%
-th row of the first matrix. Then our matrices under consideration become: 
\begin{equation*}
\left[ 
\begin{array}{cccc}
&  & \overset{k+q}{|} &  \\ 
\leftarrow & { k+q-s\rightarrow} & | &  \\ 
&  & | &  \\ 
\leftarrow & q & | & \overset{t-q}{\rightarrow} \\ 
&  & | & 
\end{array}
\right] -\left[ 
\begin{array}{cccc}
&  & \overset{k+q}{|} &  \\ 
\leftarrow & { k+q-s} & | & \overset{t-q}{\rightarrow} \\ 
&  & | &  \\ 
\leftarrow & q\rightarrow & | &  \\ 
&  & | & 
\end{array}
\right]
\end{equation*}
Now the assertion follows because there is no other choice for the first
matrix of this kind. For the general case, let the non-zero elements in the
last column of the second matrix be $l<t-q$. Then the situation is as
follows: 
\begin{equation*}
\left[ 
\begin{array}{cccc}
&  & \overset{k+q}{|} &  \\ 
\leftarrow & { k+q-s\rightarrow} & | &  \\ 
&  & | &  \\ 
\leftarrow & t-l & | & \overset{l}{\rightarrow} \\ 
&  & | & 
\end{array}
\right] -\left[ 
\begin{array}{cccc}
&  & \overset{k+q}{|} &  \\ 
\leftarrow & { k+t-s-l} & | & \overset{l}{\rightarrow} \\ 
&  & | &  \\ 
\leftarrow & q\rightarrow & | &  \\ 
&  & | & 
\end{array}
\right]
\end{equation*}
Hence we have to consider the following $(k+q)x(k+q)$ matrices: 
\begin{equation*}
\left[ 
\begin{array}{ccc}
&  &  \\ 
\leftarrow & { k+q-s} & \rightarrow \\ 
&  &  \\ 
\leftarrow & t-l & \rightarrow \\ 
&  & 
\end{array}
\right] -\left[ 
\begin{array}{ccc}
&  &  \\ 
\leftarrow & { k+t-s-l} & \rightarrow \\ 
&  &  \\ 
\leftarrow & q & \rightarrow \\ 
&  & 
\end{array}
\right]
\end{equation*}
Here the $s$-th column of the second matrix and the $k$-th column of the
first one have been raised to the power $p^{t-q-l}$. Because the exponents
are of the right form the assertion follows.
\end{proof}

For the rest of this section we recall the ring of invariants $%
(E(x_{1},...,x_{n})\otimes P[y_{1},...,y_{n}])^{GL_{n}}$ from \cite{Mui}.
Here $p>2$.

\begin{theorem}
{\rm \cite{Mui}\qua1)}\qua The algebra $(E(n)\otimes S[n])^{B_{n}}$ is a tensor product
between the polynomial algebra $B[n]$ and the $\mathbb Z /p \mathbb Z$ -module spanned by the
set of elements consisting of the following monomials: 
\begin{equation*}
M_{s;s_{1},...,s_{m}}L_{s}^{p-2};\ 1\leq m\leq n,\ m\leq s\leq n,\ \text{and 
}0\leq s_{1}<\dots<s_{m}=s-1.
\end{equation*}
Its algebra structure is determined by the following relations:\newline
{\rm a)}\qua $(M_{s;s_{1}}L_{s}^{p-2})^{2}=0$, for $1\leq s\leq n,0\leq s_{1}\leq s-1$.%
\newline
{\rm b)}\qua $M_{s;s_{1},...,s_{m}}L_{s}^{p-2}(L_{s}^{p-1})^{m-1}=\newline
(-1)^{m(m-1)/2}\prod_{q=1}^{m}(\sum
\limits_{r=s_{q}+1}^{s}M_{r;r-1}L_{r}^{p-2}h_{r+1}\dots h_{s}d_{r-1,s_{q}})$%
\hfil\break Here $1\leq m\leq n$, $m\leq s\leq n$, and $0\leq
s_{1}<\dots<s_{m}=s-1$.

{\rm 2)}\qua The algebra $(E(n)\otimes S[n])^{GL_{n}}$ is a tensor product between the
polynomial algebra $D[n]$ and the $\mathbb Z /p \mathbb Z$ -module spanned by the set of
elements consisting of the following monomials: 
\begin{equation*}
M_{n;s_{1},...,s_{m}}L_{n}^{p-2};\ 1\leq m\leq n,\ \text{and }0\leq
s_{1}<\dots<s_{m}\leq n-1.
\end{equation*}
Its algebra structure is determined by the following relations:\newline
{\rm a)}\qua $(M_{n;s_{1},...,s_{m}}L_{n}^{p-2})^{2}=0$ for$\ 1\leq m\leq n,\ $and $%
0\leq s_{1}<\dots<s_{m}\leq n-1$.\newline
{\rm b)}\qua $%
M_{n;s_{1},...,s_{m}}L_{n}^{(p-2)}d_{n,n-1}^{m-1}=(-1)^{m(m-1)/2}M_{n;s_{1}}L_{n}^{p-2}\dots M_{n;s_{m}}L_{n}^{p-2} 
$.\newline
Here $1\leq m\leq n$, and $0\leq s_{1}<\dots<s_{m}\leq n-1$.
\end{theorem}

The elements $M_{n;s_{1},...,s_{m}}$ above have been defined by Mui in \cite
{Mui}. Their degrees  are $|M_{n;s_{1},...,s_{m}}|=m+2((1+\dots
+p^{n-1})-(p^{s_{1}}+\dots +p^{s_{m}}))$ and $|L_{n}^{p-2}|=2(p-2)(1+\dots
+p^{n-1})$.

\begin{definition}
Let $S(E(n))^{B_{n}}$ be the subspace of $(E(n)\otimes S[n])^{B_{n}}$
generated by:\newline
i)\qua $M_{s;s-1}(L_{s})^{p-2}$ for $1\leq s\leq n$,\newline
ii)\qua $\prod \limits_{t=1}^{\ell}\left(
M_{s_{2t-1}+1;s_{2t-1}}(L_{s_{2t-1}+1})^{p-2}M_{s_{2t}+1;s_{2t}}(L_{s_{2t}+1})^{p-2}\right) /d_{s_{2t-1}+1,0} 
$\newline
for $0\leq s_{1}<...<s_{2\ell}\leq n-1$,\newline
iii)\qua
$M_{s_{1}+1;s_{1}}(L_{s})^{p-2}$\newline
\hbox{}\hglue 1cm$\prod \limits_{t=1}^{\ell}\left(
M_{s_{2t}+1;s_{2t}}(L_{s_{2t}+1})^{p-2}M_{s_{2t+1}+1;s_{2t+1}}(L_{s_{2t+1}+1})^{p-2}\right) /d_{s_{2t}+1,0} 
$\newline for $0\leq s_{1}<...<s_{2\ell+1}\leq n$;

and $S(E(n))^{GL_{n}}$ be the subspace of $(E(n)\otimes S[n])^{GL_{n}}$
generated by:

$M_{n;s}(L_{n})^{p-2}$ for $0\leq s\leq n-1$ ,\newline
$\prod \limits_{t=1}^{\ell}M_{n;s_{2t-1},s_{2t}}(L_{n})^{p-2}$ for $0\leq
s_{1}<...<s_{2\ell}\leq n-1$,\newline
$M_{n;s_{1}-1}(L_{n})^{p-2}\prod
\limits_{t=1}^{\ell}M_{n;s_{2t},s_{2t+1}}(L_{n})^{p-2}$ for $0\leq
s_{1}<...<s_{2\ell+1}<n$.
\end{definition}

The following lemmata provide the decomposition of $M_{n;s,m}(L_{n})^{p-2}$
in $S(E(n))^{B_{n}}\otimes B[n]$ and relations between them.

\begin{lemma}
Let $s<\ell$, then $M_{s;s-1}L_{s}^{p-2}M_{\ell;\ell-1}L_{\ell}^{p-2}$ can
be written with respect to basis elements of $B[k]\otimes S(E_{k})^{B_{k}}$.
\end{lemma}

\begin{lemma}
Let $m<s-1$, then $M_{s;\ell ,s-1}L_{s}^{p-2}M_{m;m-1}L_{m}^{p-2}$ can be
written with respect to basis elements of $S(E(k))^{B_{k}}\otimes B[k]$.
\end{lemma}

\begin{lemma}
\label{Decompo}$M_{n;s,m}(L_{n})^{p-2}=$\newline
$\sum\limits_{\substack{ s\leq q<t \\ m\leq t\leq n-1}}%
M_{q+1;q}(L_{q+1})^{p-2}M_{t+1;t}(L_{t+1})^{p-2}h_{t+2}...h_{n}(d_{q,s}d_{t,m}-d_{q_{,m}}d_{t,s})/d_{q+1,0}
$. \newline
Here $d_{i,i}=1$ and $d_{i,j}=0$ if $i<j$.
\end{lemma}

\begin{corollary}
\label{kapa}Let $\kappa =[\frac{n+1}{2}]$ and $\varepsilon =(\epsilon
_{1},...,\epsilon _{n})\in (Z/2Z)^{n}$, then $S(E(n))^{GL_{n}}$ is spanned
by at most $\kappa $ monomials: 
\begin{equation*}
M^{\varepsilon }:=\left\{ 
\begin{array}{c}
M_{n;s_{1},s_{2}}^{[\frac{\epsilon _{1}+\epsilon _{2}}{2}%
]}L_{n}^{p-2}...M_{n;s_{n-1},s_{n}}^{[\frac{\epsilon _{n-1}+\epsilon _{n}}{2}%
]}L_{n}^{p-2}\text{, if }n\text{ is even} \\ 
M_{n;s_{1}}^{\epsilon _{1}}L_{n}^{p-2}M_{n;s_{2},s_{3}}^{[\frac{\epsilon
_{2}+\epsilon _{3}}{2}]}L_{n}^{p-2}...M_{n;s_{n-1},s_{n}}^{[\frac{\epsilon
_{n-1}+\epsilon _{n}}{2}]}L_{n}^{p-2}\text{, if }n\text{ is odd}
\end{array}
\right. 
\end{equation*}
\end{corollary}

The analogue corollary holds for $S(E(n))^{B_{n}}$.

The Steenrod algebra acts naturally on $S(E(n))^{GL_{n}}\otimes D[n]$ and $%
S(E(n))^{B_{n}}\otimes B[n]$.

Let $\hat{\imath}:S(E(n))^{GL_{n}}\otimes D[n]\hookrightarrow
S(E(n))^{B_{n}}\otimes B[n]$ be the inclusion, then $\hat{\imath}%
(d^{m}M^{\varepsilon})$ means the decomposition of $d^{m}M^{\varepsilon}$ in 
$S(E(n))^{B_{n}}\otimes B[n]$.\label{inclus}

\begin{lemma}
\label{NumBock}Let $0\leq s_{1}(s_{1}^{\prime})<k_{1}(k_{1}^{\prime
})<...<s_{l^{\prime}}(s_{l^{\prime}}^{\prime})<k_{l^{\prime}}(k_{l^{%
\prime}}^{\prime})\leq n-1$. If $\sum
\limits_{0}^{n-1}m_{i}(p^{n}-p^{i})+\sum
\limits_{1}^{l^{\prime}}(p^{n}-p^{s_{i}}-p^{k_{i}})=\sum
\limits_{0}^{n-1}m_{i}^{\prime}(p^{n}-p^{i})+\sum
\limits_{1}^{l^{\prime}}(p^{n}-p^{s_{i}^{\prime}}-p^{k_{i}^{\prime}})$, then 
$s_{i}=s_{i}^{\prime}$ and $k_{i}=k_{i}^{\prime}$. Moreover, if in addition $%
0\leq k_{0}(k_{0}^{\prime})<s_{1}(s_{1}^{\prime})$ and $\sum
\limits_{0}^{n-1}m_{i}(p^{n}-p^{i})+(p^{n}-p^{k_{0}})+\sum
\limits_{1}^{l^{\prime}}(p^{n}-p^{s_{i}}-p^{k_{i}})=$ $\sum
\limits_{0}^{n-1}m_{i}^{\prime}(p^{n}-p^{i})+(p^{n}-p^{k_{0}^{\prime}})+%
\sum \limits_{1}^{l^{\prime}}(p^{n}-p^{s_{i}^{\prime}}-p^{k_{i}^{\prime}})$%
, then $s_{i}=s_{i}^{\prime}$ and $k_{j}=k_{j}^{\prime}$.
\end{lemma}

\section{Calculating the hom-duals and Adem relations}

We start this section by recalling the description of $R[n]^{\ast}$\ as an
algebra, for $p$ odd please see May \cite{C-L-M} Theorem 3.7 page 29. The
analogue Theorem for $p=2$ was given by Madsen who expressed the connection
between $R[n]^{\ast}$ \ and Dickson invariants back in 1975, \cite{Mad1}.

For convenience we shall write $I$ instead of $(I,\varepsilon)$.

Let $I_{n,i}=(\underset{i}{\underbrace{0,...,0}},\underset{n-i}{\underbrace{%
1,...,1}})$. Here $0\leq i\leq n-1$ and $n-i$ denotes the number of p-th
powers. The degree $|Q_{I_{n,i}}|=2p^{i}(p^{n-i}-1)$ [$2^{n}-2^{i}$] and the 
$exc(Q_{I_{n,i}})=0,\ $if $\ i<n$, and $1$ if $i=0$.

Let $J_{n;i}=(\underset{i}{\underbrace{\frac{1}{2},...,\frac{1}{2}}},%
\underset{n-i}{\underbrace{1,...,1}})x(\underset{i+1}{\underbrace{0,...,0,1}}%
,\underset{n-i-1}{\underbrace{0,...,0}})$. Here $\varepsilon =(\underset{i}{%
\underbrace{0,...,0}},\underset{n-i-1}{1,\underbrace{0,...,0}})$ and $0\leq
i\leq n-1$. The degree $|Q_{J_{n;i}}|=2p^{i}(p^{n-i}-1)-1$ and the $%
exc(Q_{J_{n;i}})=1$.

Let $K_{n;s,i}=(\underset{s}{\underbrace{0,...,0}},\underset{i-s}{%
\underbrace{\frac{1}{2},...,\frac{1}{2}}},\underset{n-i}{\underbrace{1,...,1}%
})x(\underset{s}{\underbrace{0,...,0}},\underset{i-s+1}{\underbrace{%
1,0,...,0,1}},\underset{n-i}{\underbrace{0,...,0}})$. Here $\varepsilon =(%
\underset{s}{\underbrace{0,...,0}},\underset{i-s}{\underbrace{1,0,...,0,1}},%
\underset{n-i}{\underbrace{0,...,0}})$ and $0\leq s<i\leq n-1$. There are
two Bockstein operations in this element: at the $s$-th and $i$-th position
from the left. The degree $|Q_{K_{n;s,i}}|=2(p^{i}(p^{n-i}-1)-p^{s})$ and
the $exc(Q_{K_{n;s,i}})=0$.

Let $O_{n,i}=(0,...,0,1,0,...,0)$, where there are $n-i$ zeros. Its degree
is $|e_{O_{n,i}}|=2p^{i-1}(p-1)$ [$2^{i-1}$] and $exc(e_{O_{n,i}})=0$. Here $%
1\leq i\leq n$.

Let $J_{n,i;i-1}=(\underset{i}{\underbrace{\frac{1}{2},...,\frac{1}{2}}},%
\underset{n-i}{\underbrace{0,...,0}})x(\underset{i}{\underbrace{0,...,0,1}},%
\underset{n-i}{\underbrace{0,...,0}})$. Here $\varepsilon =(\underset{i-1}{%
\underbrace{0,...,0}},\underset{n-i}{1,\underbrace{0,...,0}})$ and $1\leq
i\leq n$. Its degree $|Q_{J_{n,i;i-1}}|=2p^{i-1}(p-1)-1$ and the $%
exc(Q_{J_{n,i;i-1}})=1$.

Let $K_{n,i;s,i-1}=(\underset{s}{\underbrace{0,...,0}},\underset{i-s}{%
\underbrace{\frac{1}{2},...,\frac{1}{2},1}},\underset{n-i}{\underbrace{%
0,...,0}})x(\underset{s}{\underbrace{0,...,0}},\underset{i-s}{\underbrace{%
1,0,...,0,1}},\underset{n-i}{\underbrace{0,...,0}})$. Here $\varepsilon =(%
\underset{s}{\underbrace{0,...,0}},\underset{i-s}{\underbrace{1,0,...,0,1}},%
\underset{n-i}{\underbrace{0,...,0}})$ and $0\leq s<i-1\leq n-1$. Its degree 
$|Q_{K_{n,i;s,i-1}}|=2(p^{i}-p^{s}-p^{i-1})$ and the $%
exc(Q_{K_{n,i;s,i-1}})=0$.

Let $\xi_{n,0}=((Q_{0})^{n})^{\ast}=((Q^{0})^{n})^{\ast}$;\newline
$\xi
_{n,i}=(Q_{I_{n,i}})^{\ast}=(Q^{(p^{i-1}(p^{n-i}-1),%
\dots,(p^{n-i}-1),p^{n-i-1},\dots,p,1)})^{\ast}$, $0\leq i\leq n-1$;\newline
$\tau
_{n;i}=(Q_{J_{n;i}})^{\ast}=(Q^{(p^{i-1}(p^{n-i}-1),%
\dots,(p^{n-i}-1),p^{n-i-1},\dots,p,1,)x\varepsilon})^{\ast}$, $0\leq i\leq
n-1 $;\newline
$\sigma_{n;s,i}=(Q_{K_{n;s,i}})^{\ast}{=}(Q^{(p^{i-1}(p^{n-i}-1)-p^{s-1},%
\dots,p^{i-s-1}(p^{n-i}-1),\dots,\
p^{n-i}-1,p^{n-i-1},\dots,p,1)x\varepsilon})^{\ast}$, $0\leq s<i\leq n-1$;%
\newline
$\zeta _{n,i}=\left( e_{_{^{O_{n,i}}}}\right)
^{\ast}=(e^{(p^{i-2}(p-1),\dots ,(p-1),1,0,\dots,0)})^{\ast}$, $1\leq i\leq
n $;\newline
$\nu_{n,i;i-1}=(e_{J_{n,i;i-1}})^{\ast}=(e^{(p^{i-2}(p-1),\dots,(p-1),1,0,%
\dots ,0)x\varepsilon})^{\ast}$, $1\leq i\leq n$;\newline
$\upsilon_{n,i;s,i-1}=(e_{K_{n,i;s,i-1}})^{\ast}=$ \newline
$(e^{(p^{i-1}(p^{n-i}-1)-p^{s-1},\dots,p^{i-s-1}(p^{n-i}-1),%
\dots,p^{n-i}-1,p^{n-i-1},\dots,p,1,0,\dots ,0)x\varepsilon})^{\ast}$,
\newline $0\leq s<i-1\leq n-1$.

\begin{theorem}[Madsen $p=2$, May $p>2$]
As an $A$ algebra $R[n]^{\ast}\cong$ free associative commutative algebra
generated by $\{\xi_{n,i},\tau_{n;i}$, and $\sigma_{n;s,i}\ |\ 0\leq i\leq
n-1$, and $0\leq s<i\}$, [$\{\xi _{n,i}\ |\ 0\leq i\leq n-1\}$], modulo the
following relations:\newline
{\rm a)}\qua $\tau_{n;i}\ \tau_{n;i}=0$. \newline
{\rm b)}\qua $\tau_{n;s}\tau_{n;i}=\sigma _{n;s,i}\xi_{n,0}$. Here $0\leq s<i\leq n-1$%
. \newline
{\rm c)}\qua $\tau_{n;s}\tau_{n;i}\tau_{n;j}=\tau_{n;s}\sigma_{n;i,j}\xi_{n,0}.$ Here $%
0\leq s<i<j\leq n-1$.\newline
{\rm d)}\qua $\tau_{n;s}\tau_{n;i}\tau_{n;j}\tau_{n;k}=\sigma
_{n;s,i}\sigma_{n;j,k}\xi_{n,0}^{2}.$ Here $0\leq s<i<j<k\leq n-1$.
\end{theorem}

\begin{theorem}
{\rm\cite{Kech}}\qua\label{Phi}$R[n]^{\ast }\approxeq S(E(n))^{GL_{n}}\otimes D[n]$ [%
$R[n]^{\ast }\approxeq D[n]$] and $T[n]^{\ast }\approxeq
S(E(n))^{B_{n}}\otimes B[n]$ [$T[n]^{\ast }\approxeq B[n]$] as algebras over
the Steenrod algebra {and the isomorphism }$\Phi $ is given by $\Phi ${$(\xi
_{n,i}=(Q_{I_{n,i}})^{\ast })=d_{n,n-i}$, }$\Phi $$(\tau
_{n;i}=(Q_{J_{n;i}})^{\ast })=M_{n;i}(L_{n})^{p-2}$, $\Phi ${$($}$\sigma
_{n;s,i}=(Q_{K_{n;s,i}})^{\ast })=M_{n;s,i}(L_{n})^{p-2}$. Here $0\leq i\leq
n-1$ and $0\leq s<i$.\newline
$\Phi $$(\zeta _{n,i}=\left( e_{^{O_{n,i}}}\right) ^{\ast })=h_{i}$, $\Phi ${%
$($}$\nu _{n,i;i-1}=(e_{J_{n,i;i-1}})^{\ast })=M_{i;i-1}(L_{i})^{p-2}$, $%
\Phi $$(\upsilon _{n,i;s,i-1}=(e_{K_{n,i;s,i-1}})^{\ast
})=(M_{s+1;s}(L_{s+1})^{p-2}M_{i;i-1}(L_{i})^{p-2})/d_{s+1,0}$.{\ Here $%
1\leq i\leq n$ and }$0\leq s<i-1$.
\end{theorem}

Under isomorphism $\Phi$ in Theorem \ref{Phi} we identify $R[n]^{\ast}$ with 
$S(E(n))^{GL_{n}}\otimes D[n]$ and $B[n]^{\ast}$ with $S(E(n))^{B_{n}}%
\otimes B[n]$.

The set \c{T}$[n]$ and \c{R}$[n]$ of admissible monomials in $T[n]$ and $%
R[n] $ provide vector space bases respectively. Let $\theta:\c{R}[n%
]\rightarrowtail $\c{T}$[n]$ be the map given by 
\begin{equation*}
\theta(Q_{I})=e_{I}
\end{equation*}
The image of the dual of these bases are denoted by \c{T}$[n]^{\ast}$ in $%
\Phi(T[n])^{\ast}=S(E(n))^{B_{n}}\otimes B[n]$ and \c{R}$[n]^{\ast}$ in $%
\Phi(R[n])^{\ast}=S(E(n))^{GL_{n}}\otimes D[n]$. Of course there are also
the bases of monomials which are denoted by \ss$_{n}(S(E(n))^{B_{n}}\otimes
B[n])$ and \ss$_{n}(S(E(n))^{GL_{n}}\otimes D[n])$ respectively. We shall
note that \c{T}$[n]^{\ast}=$\ss$_{n}(S(E(n))^{B_{n}}\otimes B[n])$.

The decomposition relations between the other two bases are not obvious and
this is the first topic of this section. Campbell, Peterson and Selick
provided a method to pass from \ss $_{n}(D[n])$ to \c{R}$[n]^{\ast }$ in 
\cite{C-P-S}. We shall establish some machinery to work with those bases.

\begin{definition}
Let $\chi _{\min }$ and $\chi _{\max }$ be the set functions from \ss $%
(S(E(n))^{GL_{n}}\otimes D[n])$ (\ss $(B[n]\otimes S(E_{n})^{B_{n}})$ ) to
the monoid $\langle \mathbb N ,\frac{1}{2}\rangle ^{n}\times (\mathbb Z /2\mathbb Z )^{n}$ given by\newline
1)\qua $\chi _{\min }(d_{n,i})=I_{n,i}$, $\chi _{\max
}(d_{n,i})=(p^{n-i},...,p^{n-i},0,...,0)x(0,...,0)$;\newline
2)\qua  $\chi _{\min }(M_{n;s}L_{n}^{(p-2)})=J_{n;s}$,\newline 
$\chi _{\max
}(M_{n;s}L_{n}^{(p-2)})=(\underset{s}{\underbrace{\frac{1}{2},...,\frac{1}{2}%
}},\underset{n-s-1}{\underbrace{1\frac{1}{2},...,1\frac{1}{2}}}%
,1)x(0,...,0,1)$;\newline
3)\qua  $\chi _{\min }(M_{n;s,m}L_{n}^{(p-2)})=K_{n;s,m}$ and \newline
$\chi _{\max }(M_{n;s,m}L_{n}^{(p-2)})=(\underset{m}{\underbrace{0,...,0}},%
\underset{n-m-1}{\underbrace{1\frac{1}{2},...,1\frac{1}{2}}},1)x(\underset{m%
}{\underbrace{0,...,0}},\underset{n-m}{\underbrace{1,0,...,0,1}})$; \newline
and the rule $\chi _{\min }(dd^{\prime }MM^{\prime })=\chi _{\min }(d)+\chi
_{\min }(d^{\prime })+\chi _{\min }(M)+\chi _{\min }(M^{\prime })$. Here $d$%
, $d^{\prime }\in $\ss $(D[n])$ and $M$, $M^{\prime }\in $\ss $%
(S(E(n))^{GL_{n}})$. The same holds for $\chi _{\max }$.
\end{definition}

Note that the function $\chi_{\min}$ always provides an admissible element
and $\hat{\imath}(d_{n,i})$ contains a monomial with a unique admissible
sequence, namely $h^{\chi_{\min}(d_{n,i})}$, and a monomial with a unique
maximal sequence, namely $h^{\chi_{\max}(d_{n,i})}$. The same is true for
elements $M_{n;s-1}L_{n}^{p-2}$ and $M_{n;s,m}L_{n}^{p-2}$. Moreover, $\hat{%
\imath }(d_{n}^{m}M)$ might contain a number of monomials with admissible
sequences and this is the main point of investigation because of its
applications in \cite{kech4}. Namely, those monomials provide possible
candidates for $(d_{n}^{m}M)^{\ast }$. Primitives in $R$ are well known and
so are their duals as generators in $R^{\ast }$. But it is not the case for
their expression with respect to the Dickson algebra. On the other hand, the
action on the Dickson algebra is well known on $S(E(n))^{GL_{n}}\otimes D[n]$
and hence it is easier to compute the annihilator ideal in the mod-$p$
cohomology of a certain finite loop space.

\begin{definition}
Let $\Psi $ be the correspondence between \ss $_{n}(S(E(n))^{GL_{n}}\otimes
D[n])$ and \c{R}$[n]$ given by $d\longmapsto \Psi (d)=Q_{\chi _{\min }(d)}$
and the corresponding one between \ss $_{n}(S(E(n))^{B_{n}}\otimes B[n])$
and \c{T}$[n]$ denoted by $\Psi _{T}$ where\newline
$\Psi _{T}(h^{J}M^{\varepsilon })=e_{J+\epsilon
_{1}J_{n;s_{1}}+\sum\limits_{t}[\frac{\epsilon _{2t}+\epsilon _{2t+1}}{2}%
]K_{n;s_{2t},s_{2t=1}}}$. 
\end{definition}

The maps $\Psi $ and $\Psi _{T}$ are set bijections.

Let $\iota$ be the map 
\begin{equation}
\iota:\text{\ss}_{n}(S(E(n))^{GL_{n}}\otimes D[n])\rightarrowtail\text{\ss }%
_{n}(S(E(n))^{B_{n}}\otimes B[n])  \label{kkk}
\end{equation}
defined by $\iota(d)=h^{\chi_{\min}(d)}$.

Note that $e_{\chi_{\max}(d)},e_{\chi_{\min}(d)}\in$\c{T}$[n]$. The
following diagram is commutative. 
{\small\begin{equation*}
\begin{array}{ccccc}
\text{\ss}_{n}(S(E(n))^{B_{n}}\otimes B[n]) & \overset{\iota}{\leftarrowtail}
& \text{\ss}_{n}(S(E(n))^{GL_{n}}\otimes D[n]) & \overset{\chi_{\min}}{%
\rightarrow} & \langle \mathbb N ,\frac{1}{2}\rangle ^{n}\times(\mathbb Z /2\mathbb Z )^{n} \\ 
\begin{array}{cc}
\;\downarrow & \Psi_{T}
\end{array}
&  & 
\begin{array}{cc}
\;\downarrow & \Psi
\end{array}
& \swarrow &  \\ 
\c{T}[n] & \overset{\theta}{\leftarrowtail} & \c{R}[n] &  & 
\end{array}
\end{equation*}}

\begin{definition}
A monomial in $\ss_{n}(S(E(n))^{B_{n}}\otimes B[n])$ is called admissible if
it is an element of $\iota\left( \text{\ss}_{n}(S(E(n))^{GL_{n}}\otimes
D[n])\right) $.
\end{definition}

\begin{lemma}
Let $h^{J}M^{\varepsilon }\in S(E(n))^{B_{n}}\otimes B[n]$. The following
are equivalent:\newline
{\rm i)}\qua $h^{J}M^{\varepsilon }$ is admissible;\newline
{\rm ii)}\qua $j_{t}\leq j_{t+1}$ for $t=1,...,n-1$ and $h^{J}$ is divisible by $%
\prod\limits_{t=0}^{l}(h_{s_{2t+1}+2}...h_{n})^{\epsilon _{2t+1}}$ for $%
\kappa $ odd (see \ref{kapa}); or $\prod%
\limits_{t=1}^{l}(h_{s_{2t}+2}...h_{n})^{\epsilon _{2t}}$, otherwise. If $%
s_{2t+1}+2$ or $s_{2t+2}+2=n+1$, then the corresponding product must be $1$.

{\rm iii)}\qua $\rho(\Psi_{T}(h^{J}M^{\varepsilon}))$ is admissible in $R[n]$.
\end{lemma}

\begin{proof}
This follows from the following relation:\newline
$M_{k;s}L_{k}^{p-2}=M_{s+1;s}L_{s+1}^{p-2}h_{s+1}...h_{k}+\sum%
\limits_{t=2}^{k-s}M_{s+t;s+t-1}L_{s+t}^{p-2}d_{s+t-1,s}h_{s+t+1}...h_{k}$.
Explicitly, if $h^{J^{\prime
}}=h^{J}/\prod\limits_{t=0}^{l}(h_{s_{2t+1}+2}...h_{n})^{\epsilon _{2t+1}}$%
, then $\chi _{\min }(d^{m}M^{\varepsilon })=(J^{\prime },\varepsilon )$.
Here\newline
$d^{m}M^{\varepsilon
}=\prod\limits_{i=0}^{n-1}d_{n,i}^{m_{i}}M_{n;s_{1}}^{\epsilon
_{1}}L_{n}^{p-2}M_{n;s_{2},s_{3}}^{[\frac{\epsilon _{2}+\epsilon _{3}}{2}%
]}L_{n}^{p-2}...M_{n;s_{n-1},s_{n}}^{[\frac{\epsilon _{n-1}+\epsilon _{n}}{2}%
]}L_{n}^{p-2}$ and $m_{t}=j_{t}^{\prime }-j_{t-1}^{\prime }$, $%
m_{0}=j_{0}^{\prime }$.
\end{proof}

Firstly, we shall show that $\rho^{\ast}\equiv\hat{\imath}$, $\hat{\imath}$
as in \ref{inclus}, i.e.\ for any $e_{I}\in T[n]$ and $d^{m}M^{\varepsilon}$, 
\begin{equation*}
\langle d^{m}M^{\varepsilon},\rho(e_{I})\rangle =\langle \hat{\imath}(d^{m}M^{%
\varepsilon}),e_{I}\rangle .
\end{equation*}
Here, $\langle -,-\rangle $ is the Kronecker product. This is done by studying all
monomials in $T[n]$ which map to primitives in $R[n]$ after applying Adem
relations. \label{AAdem}

Let $n(mx\varepsilon )=\sum m_{i}+\kappa $. Let $\psi _{n(mx\varepsilon
)}:R[n]\rightarrow \overset{n(mx\varepsilon )}{\otimes }R[n]$ be the
iterated coproduct $n(mx\varepsilon )$ times. Let $J$ be admissible, $\rho
e_{J}=Q_{J}$, then 
\begin{align*}
\psi Q_{J}& =\psi \rho e_{J}=\rho \psi e_{J}=\rho (\Sigma \pm
e_{J_{1}}\otimes \cdots \otimes e_{J_{n(mx\varepsilon )}}),\quad \Sigma
J_{i}=J \\
\psi _{n(mx\varepsilon )}Q_{J}& =\Sigma a_{J_{1},...,J_{n(mx\varepsilon
)}}Q_{J_{1}^{\prime }}\otimes \cdots \otimes Q_{J_{n(mx\varepsilon
)}^{\prime }}.
\end{align*}
Since $J_{i}$ may not be in admissible form, after applying Adem relations
we have $J_{i}^{\prime }\leq J_{i}$.\newline
$\langle d^{m}M^{\varepsilon },\rho
e_{I}\rangle =\langle \prod\limits_{i=0}^{n-1}d_{n,i}^{m_{i}}M^{\varepsilon },\psi
_{n(mx\varepsilon )}\rho
e_{I}\rangle =\langle \prod\limits_{i=0}^{n-1}d_{n,i}^{m_{i}}M^{\varepsilon },\rho \psi
_{n(mx\varepsilon )}e_{I}\rangle =$\newline
$\langle \prod\limits_{i=0}^{n-1}d_{n,i}^{m_{i}}M^{\varepsilon
},\sum\limits_{I_{j}}\bigotimes\limits_{j}^{n(mx\varepsilon )}\rho
e_{I_{j}}\rangle $ $=\sum\limits_{I_{j}}\prod\limits_{j}^{n(m)}\langle d_{n,i}^{j},\rho
e_{I_{j}}\rangle \prod\limits_{j}^{n(\varepsilon )}\langle M_{n;s_{j-1},s_{j}}^{[\frac{%
\epsilon _{j-1}+\epsilon _{j}}{2}]}L_{n}^{p-2},\rho e_{I_{j}}\rangle $.

\begin{lemma}
Let $d^{m}=\prod\limits_{i=1}^{n}d_{n,i}^{m_{i}}$. Then $\iota \left(
d^{m}\right) =\prod\limits_{t=1}^{n}h_{t}^{\sum\limits_{i=0}^{t-1}m_{i}}$ and%
\newline
$\left( \iota \left( d^{m}\right) \right) ^{\ast
}=e_{m_{0}}e_{m_{0}+m_{1}}...e_{m_{0}+...m_{n-1}}$.\label{min-seq}
\end{lemma}

\begin{lemma}
Let $I=\chi_{\max}(d_{n,n-i})$, then $\rho(e_{I})=Q_{I_{n,n-i}}=\Psi\left(
d_{n,n-i}\right) $ in $R[n]$.
\end{lemma}

\begin{proof}
By direct computation.
\end{proof}

\begin{lemma}
\label{A}Let $e_{I}\in T[n]$ be such that $e_{I}=\Phi^{-1}\left(
\prod\limits_{s=1}^{i}h_{j_{s}}^{p^{n-i+s-j_{s}}}\right) $ in (\ref{formula}%
). Here $1\leq j_{1}<...<j_{i}\leq n$. Then $\rho(e_{I})=Q_{n,n-i}=\Psi%
\left( d_{n,n-i}\right) $ in $R[n]$.
\end{lemma}

\begin{proof}
The sequence $I$ is given by: 
\begin{equation*}
\left( \underset{j_{1}}{\underbrace{0,\cdots,0,p^{n-i+1-j_{1}}}},\underset{%
j_{i-2}-j_{1}}{\underbrace{0,\cdots,0,p^{n-i+2-j_{2}}}},\cdots,\underset{%
j_{i}-j_{i-1}}{\underbrace{0,\cdots,0,p^{n-j_{i}}}},\underset{n-j_{i}}{%
\underbrace{0,\cdots,0}}\right)
\end{equation*}
Please note the analogy between $I$ above and the corresponding row of a
matrix in $A_{n,n-i}$ in section 3. Here $p^{m}:=0$, whenever $m<0$. We
shall work out the first steps to describe the idea of the proof. First, we
consider the last $n-i+1$ elements of $\chi_{\max}(d_{n,n-i})$: $%
(p^{n-i},0,...,0)$ which becomes $(\underset{j_{i}-j_{i-1}}{\underbrace{%
0,\cdots,0,p^{n-j_{i}}}},\underset{n-j_{i}}{\underbrace{0,\cdots,0}})$. Thus
applying Adem relations on certain positions on $Q_{\chi_{\max}(d_{n,n-i})}$%
, $Q_{I}$ is obtained and the lemma follows.
\end{proof}

\begin{proposition}
\label{B}Let $e_{I}\in T[n]$ be the $\hom$-dual of a monomial $h^{J}\in T[n]$
such that $|h^{J}|=2\left( p^{n}-p^{n-i}\right) $ and $h^{J}$ is not a
summand in (\ref{formula}). Then $\rho(e_{I})=0$ in $R[n]$.
\end{proposition}

\proof
Please see:\newline 
\url{http://www.maths.warwick.ac.uk/agt/ftp/aux/agt-4-13/full.ps.gz}\hfil$\sq$
\medskip

Now the following theorem is easily deduced because $R[n]$ is a coalgebra,
the map $\rho$ is a coalgebra map, and primitives which do not involve
Bockstein operations have been considered.

\begin{theorem}
\label{A'dem}Let $\rho^{\prime}$ be the restriction of $\rho$ \ between the
subcoalgebras $T^{\prime}[n]$ and $R^{\prime}[n]$ where no Bockstein
operations are allowed. Let $\hat{\imath}^{\prime}:D[n]\hookrightarrow B[n]$
be the natural inclusion. Then $(\rho^{\prime})^{\ast}\equiv\hat{\imath }%
^{\prime}$, i.e.\ for any $e_{I}\in T[n]$ and $d^{m}=\prod%
\limits_{i=0}^{n-1}d_{n,i}^{m_{i}}\in D[n]$, 
\begin{equation*}
\langle d^{m},\rho^{\prime}(e_{I})\rangle =\langle \hat{\imath}^{\prime}(d^{m}),e_{I}\rangle .
\end{equation*}
\end{theorem}

We shall extend last Theorem to cases including Bockstein operations as well.

Please recall that $\chi _{\min }(M_{n;s}L_{n}^{(p-2)})=J_{n;s}$.

\begin{proposition}
{\rm a)}\qua Let $J=J_{n,t;t-1}+(I_{t-1,s}^{\prime}\oplus0_{n-t+1})+I_{n,t}$ such that 
$\rho^{\prime}(e_{I_{t-1,s}^{\prime}})=Q_{I_{t-1,s}}$ for $s+1\leq t\leq n$.
Then $\rho(e_{J})=Q_{J_{n;s}}=\Psi(M_{n;s}L_{n}^{p-2})$.

{\rm b)}\qua  Let $J$ be a sequence of length $n$ such that $|J|=2(p^{n}-p^{s})-1$ and $%
J$ is not of the form described in a), then $\rho(e_{J})=0$.

{\rm c)}\qua $\hat{\imath}(M_{n;s}L_{n}^{p-2})=\rho^{\ast}(M_{n;s}L_{n}^{p-2})$.
\end{proposition}

\proof
Please see:\newline
\url{http://www.maths.warwick.ac.uk/agt/ftp/aux/agt-4-13/full.ps.gz}\hfil$\sq$\medskip

\begin{lemma}
{\rm a)}\qua Let the sequences $K_{n,t+1;q,t}$ and $I_{n,t+1}$, then $\rho
e_{(K_{n,t+1;q,t}+I_{n,t+1})}=Q_{K_{n;q,t}}$.

{\rm b)}\qua Let the sequence $K=K_{n;q,t}+(I_{q}^{"}\oplus0_{n-q})+(I_{t}^{\prime
}\oplus0_{n-t})$ such that $I_{t}^{\prime}=I_{q}^{\prime}\oplus
I_{t-q}^{\prime}$, $\rho(e_{I_{q}^{"}})=Q_{I"_{q,s}}$ and $%
\rho(e_{I_{t}^{\prime}})=Q_{I_{t,m}}$. If we allow Adem relations everywhere
in the first $t$ positions except at positions between $q$ and $q+1$ from
left, then $\rho^{\prime}(e_{K})=e_{K^{\prime}}$ where $K^{%
\prime}=K_{n;q,t}+(I_{q,s}^{"}\oplus0_{n-q})+p^{t-q-m_{2}}(I_{q,m_{1}}^{%
\prime}\oplus0_{n-q})+(0_{q}\oplus I_{t-q,m_{2}}^{\prime}\oplus0_{n-t})$ or $%
K^{\prime}=K_{n;q,t}+p^{t-q-m_{2}}(I_{q,s+m_{1}-q}^{\prime}%
\oplus0_{n-q})+(0_{q}\oplus I_{t-q,m_{2}}^{\prime}\oplus0_{n-t})$. For the
first case $\rho(e_{I_{t}^{\prime}})=Q_{I_{t-q,m_{2}}^{\prime}}$, $%
\rho(e_{I_{q}^{\prime}})=Q_{p^{t-q-m_{2}}I_{q,m_{1}}^{\prime}}$ and $%
m=m_{1}+m_{2}$, and for the second $s+m_{1}\geq q$ and $%
\rho(e_{I"_{q}+I_{q}^{\prime}})=Q_{p^{t-q-m_{2}}I_{q,s+m_{1}-q}^{\prime}}$.
\end{lemma}

\begin{proof}
This is an application of theorem \ref{A'dem}.
\end{proof}

\begin{proposition}
{\rm a)}\qua Let $K=K_{n,t+1;s,t}+(I_{t}^{\prime}\oplus0_{n-t})+I_{n,t+1}$ such that $%
\rho^{\prime}(e_{I_{t}^{\prime}})=Q_{I_{t,m}}$ for $m\leq t\leq n-1$. Then $%
\rho(e_{K})=Q_{K_{n;s,m}}=\Psi(M_{n;s,m}L_{n}^{p-2})$.

{\rm b)}\qua Let $K=K_{n,m+1;t,m}+(I_{t}^{\prime}\oplus0_{n-t})+I_{n,m+1}$ such that $%
\rho^{\prime}(e_{I_{t}^{\prime}})=Q_{I_{t,s}}$ for $s\leq t\leq m-1$. Then $%
\rho(e_{K})=Q_{K_{n;s,m}}=\Psi(M_{n;s,m}L_{n}^{p-2})$.

{\rm c)}\qua Let $K=K_{n,t+1;q,t}+I+I_{n,t+1}$ for $m\leq q<t\leq n-1$ with $%
I=I^{\prime}+I"$, $I^{\prime}=(I_{q}^{\prime}\oplus0_{n-q})$, $%
I"=(I_{t}^{"}\oplus0_{n-t})$ such that: $\rho^{%
\prime}(e_{I_{t}^{"}})=Q_{I_{t,m}}$ and $\rho^{\prime}(e_{I_{q}^{^{%
\prime}}})=Q_{I_{q,s}}$ and not of the form $\rho^{%
\prime}(e_{I_{t}^{"}})=Q_{I_{t,s}}$ and $\rho^{\prime}(e_{I_{q}^{^{%
\prime}}})=Q_{I_{q,m}}$. Then $\rho(e_{K})=Q_{K_{n;s,m}}=\Psi
(M_{n;s,m}L_{n}^{p-2})$.

{\rm d)}\qua Let $K$ be a sequence of length $n$ such that $|K|=2(p^{n}-p^{s}-p^{m})$
and $K$ is not of the form described in a), b) and c) above, then $\rho
(e_{K})=0$.

{\rm e)}\qua $\hat{\imath}(M_{n;s,m}L_{n}^{p-2})=\rho^{\ast}(M_{n;s,m}L_{n}^{p-2})$.
\end{proposition}

\proof
Please see:\newline 
\url{http://www.maths.warwick.ac.uk/agt/ftp/aux/agt-4-13/full.ps.gz}\hfil$\sq$\medskip

\begin{theorem}
Let $\rho:T[n]\rightarrow R[n]$ be the map which imposes Adem relations. Let 
$\hat{\imath}:S(E(n))^{GL_{n}}\otimes D[n]\hookrightarrow
S(E(n))^{B_{n}}\otimes B[n]$ be the natural inclusion. Then $%
\rho^{\ast}\equiv\hat{\imath}$, i.e.\ for any $e_{I}\in T[n]$ and $%
d^{m}M^{\varepsilon}\in S(E(n))^{GL_{n}}\otimes D[n]$, 
\begin{equation*}
\langle d^{m}M^{\varepsilon},\rho (e_{I})\rangle =\langle \hat{\imath} (d^{m}M^{%
\varepsilon}),e_{I}\rangle .
\end{equation*}
\end{theorem}

\begin{theorem}
\label{dual}Let $d^{m}M^{\varepsilon }$ be an element of \ss $%
_{n}(S(E(n))^{GL_{n}}\otimes D[n])$, then the following algorithm calculates
its image in $R[n]^{\ast }$: 
\begin{equation*}
d^{m}M^{\varepsilon }=\sum\limits_{J\geq \chi _{\min
}(d^{m})}\langle d^{m},Q_{J}\rangle (Q_{(J+\chi _{\min }(M^{\varepsilon }))})^{\ast }
\end{equation*}
{\rm1)}\qua Find all elements $Q_{J}$ in $R[n]$ such that $|d^{m}|=|Q_{J}|$ and $%
J>\chi _{\min }(d^{m})$,  i.e.\ solve the Diophantine equation $%
\sum\limits_{0}^{n-1}k_{i}(p^{n}-p^{i})=|d^{m}|$ for $%
(k_{0},...,k_{n-1})>(m_{0},...,m_{n-1})$. For each such a sequence $J$, let $J(1)=J-m_0(1,...,1)$ and consider $\Psi
^{-1}(Q_{J(1)})=d^{J^{\prime }(1)}$ in $D[n]$.

{\rm2)}\qua Let $d^{m}M^{\varepsilon}=(Q_{\chi_{\min}(d^{m}M^{\varepsilon})})^{\ast}$.

{\rm3)}\qua Let $d^{m(1)}=\frac{\rm d^m}{\rm d_{n,0}^{m_0}}$ and $d^{K}$ be an element in step 1) corresponding to the biggest
sequence among those which have not been considered yet. If $%
d^{K(1)}=d^{m(1)}$, then $\alpha_{(K)}=\langle d^{m},Q_{K}\rangle =1$. Otherwise, proceed
as follows: find the coefficient, $\alpha_{(K)}$, of $\iota(d^{K(1)})$ in $%
\hat{\imath}(d^{m(1)})$, $\alpha_{(K)}=\langle d^{m},Q_{K}\rangle $. Then add $\alpha
_{(K)}(Q_{K+\chi_{\min}(M^{\varepsilon})})^{\ast}$ in $d^{m}M^{\varepsilon}$.

{\rm4)}\qua Repeat step 3).
\end{theorem}

\begin{proof}
Since $R[n]^{\ast }$ is a free module over $D[n]$ with basis all elements
which involve Bockstein operations, the computation of $d^{m}M^{\varepsilon }
$ reduces to that of $d^{m}$, i.e.\
\begin{gather*}d^{m} =\sum\limits_{J\geq \chi _{\min
}(d^{m})}\langle d^{m},Q_{J}\rangle (Q_{(J)})^{\ast }\qquad\Rightarrow\\
d^{m}M^{\varepsilon } =\sum\limits_{J\geq \chi _{\min
}(d^{m})}\langle d^{m},Q_{J}\rangle (Q_{(J+\chi _{\min }(M^{\varepsilon }))})^{\ast }\end{gather*}
Let $d^{m}=\sum \alpha_{(I)}(Q_{I})^{\ast}$ and $n(m)=\sum%
\limits_{t=0}^{n-1}m_{t}$. Because of the definition of the $\hom$-dual, we
have : $\langle d^{m},Q_{\chi_{\min}(d^{m})}\rangle =1$ and $\langle d^{m},Q_{I}\rangle =a_{(I)}\neq0$
for a sequence $I$ such that in the $n(m)$-times iterated coproduct: 
\begin{equation*}
\psi Q^{I}=\sum_{\Sigma J_{t}=I}e_{J_{1}}\otimes...\otimes e_{J_{n(m)}}%
\overset{Adem}{=}\sum a_{I_{1},...,I_{n(m)}}Q_{I_{1}}\otimes...\otimes
Q_{I_{n(m)}}
\end{equation*}
$a_{(I)}\bigotimes \limits_{t=0}^{n-1}\bigotimes
\limits_{1}^{m_{t}}(Q_{I_{n,t}})$ is a summand. Thus $I\geq\chi_{%
\min}(d^{m}) $. Let $I_{1}>\cdots>I_{l}>\chi_{\min}(d^{m})$ be all sequences
such that $|Q_{I_{t}}|=|Q_{\chi_{\min}(d^{m})}|$.

We quote from May page 20: if for each $d^{m}M^{\varepsilon}$ we associate
its coefficients $a_{(I)}$ as a matrix $(a_{\chi_{\min}(d^{m}M^{%
\varepsilon}),(I)})$, then this matrix is upper triangular with ones along
the main diagonal. This allows us to express one basis element $%
d^{m}M^{\varepsilon}$ with respect to the dual basis of admissible monomials.

We consider the first sequence $I_{1}$. Our task is to evaluate $\alpha
_{(I_{1})}$. Let $\psi Q_{I_{1}}$ be the iterated coproduct applied $n(m)$%
-times. We shall write $I_{1}$ as a sum of $n(m)$ sequences such that each
of them is a primitive element of $R[n]$ equals to one of those involved in $%
\chi_{\min}(d^{m})$. This is possible, since $n(m)\geq
n(\chi_{\min}(\Psi^{-1}(Q_{I_{1}})))$. The common element $d_{n,0}^{m_{0}}$
between $\Psi^{-1}(Q_{(I_{1})})$ and $d^{m}$ does not change the coefficient 
$\alpha_{(I_{1})}$, because no Adem relation can reduce $Q_{I_{n,0}}$ to a
smaller sequence. Instead, we consider $Q_{I_{1}-m_{0}I_{n,0}}$ ($%
d^{J_{1}}=\Psi^{-1}(Q_{(I_{1})})/d_{n,0}^{m_{0}}$) and $Q_{(\chi_{%
\min}(d^{m})-m_{0}I_{n,0})}$ ($d^{m(1)}=d^{m}/d_{n,0}^{m_{0}}$). Now the
iterated coproduct is applied $n(m(1))$-times.

For the second part of step 3), we use $\psi\rho=\rho\psi$, lemma \ref{A}
and proposition \ref{B}. All elements $e_{I}\in T[n]$, which have the
property $\rho e_{I}=Q_{I_{n,n-i}}$, are known. Moreover, the dual of those
elements, $(e_{I})^{\ast}\in B[n]$, are summands in $\hat{\imath}(d_{n,n-i}) 
$. Using commutativity in $D[n]$ induced by symmetry in coproduct, we deduce
that the required coefficient is the coefficient of $\iota(d^{J_{1}})$ in $%
\hat{\imath }(d^{m(1)})$.
\end{proof}

\begin{remark}
Suppose that $(Q_{I})^{\ast}$ is to be expressed with respect to \ss $%
_{n}(S(E(n))^{GL_{n}}\otimes D[n])$, then one starts with the biggest
sequence, say $K(1)$, $\Psi^{-1}(Q_{K(1)})=(Q_{K(1)})^{\ast}$, then
substitutes in the next element $\Psi^{-1}(Q_{K(2)})=(Q_{K(2)})^{\ast
}+a_{K(2),K(1)}(Q_{K(1)})^{\ast}$ $\Rightarrow$ $(Q_{K(2)})^{\ast}=\Psi
^{-1}(Q_{K(2)})-a_{K(2),K(1)}\Psi^{-1}(Q_{K(2)})$ and so on.
\end{remark}

Let us make some comments. If the degree $m$ of a monomial $d^{m}$ is quite
high, then there exist many elements of the same degree such that the dual
of their images under $\Phi$ do not appear in $d^{m}$ for a variety of
reasons. We shall give a refinement of the algorithm described above through
the next lemmas.

\begin{definition}
Let $d^{m}=\prod \limits_{i=0}^{n-1}d_{n,i}^{m_{i}}$ be a monomial in the
Dickson algebra and $m_{i}=\sum\limits_{t=0}^{\ell_{i}}a_{i,t}p^{t}$. Let $%
i_{0}=\max\{i\;|\;m_{i}\neq0\}$ and $0\leq t<i_{0}$. Let $\delta(t)$ be a
positive integer such that $t\leq\gamma(s)\leq n-1$ for $s=1,...,\delta(t) $
and $\sum\limits_{1}^{\gamma(\delta(t))}(n-\gamma(s))=n-t$. Let also $%
\ell_{(t,\gamma (1),...,\gamma(\delta(t)))}=\max\{\gamma(s)-(\sum
\limits_{1}^{s}\gamma(j))+(s-1)n\;|\;s=1,...,\gamma(\delta(t))\} $ and $%
0\leq c\leq\min
\{\ell_{\gamma(s)}-\ell_{(t,\gamma(1),...,\gamma(\delta(t)))}+\sum
\limits_{j=1}^{s-1}(n-\gamma(j))\}$. We define 
\begin{equation*}
\zeta(t,\gamma(1),...,\gamma(\delta(t)),c,\mu)=\left\{ 
\begin{array}{l}
\prod\limits_{s=2}^{\gamma(\delta(t))}\binom{a_{\gamma(s),c+\ell
_{(t,\gamma(1),...,\gamma(\delta(t)))}-\sum\limits_{j=1}^{s-1}(n-\gamma(j))}%
}{\mu}\text{,}\; \\ 
\text{if }0\leq m_{\gamma(1)}-\mu p^{c+\ell_{(t,\gamma(1),...,\gamma
(\delta(t)))}} \\ 
0\text{, otherwise}
\end{array}
\right. .
\end{equation*}
Here $1\leq\mu\leq\min\{a_{\gamma(s),c+\ell_{(t,\gamma(1),...,\gamma
(\delta(t)))}-\sum\limits_{j=1}^{s-1}(n-\gamma(j))}\;|\;s=2,...,\delta(t)\}$.
\end{definition}

\begin{proposition}
Let $d^{m}=\prod \limits_{i=0}^{n-1}d_{n,i}^{m_{i}}$ be a monomial in the
Dickson algebra as above. Then $d^{m}$ contains $$\left(
\Psi(d^{m}d_{n,t}^{p^{\ell_{(t,\gamma(1),...,\gamma
(\delta(t)))}}}/\prod\limits_{1}^{\gamma(\delta(t))}d_{n,\gamma(s)}^{p^{%
\ell_{(t,\gamma(1),...,\gamma(\delta(t)))}-\sum\limits_{j=1}^{s-1}(n-%
\gamma(j))}})\right) ^{\ast}$$ with coefficient 
\begin{align*}
& \sum\limits_{\gamma(1),...,\gamma(\delta(t))}\zeta(t,\gamma(1),...,\gamma
(\delta(t)),0,1)+ \\
&
\sum\limits_{\gamma^{\prime}(1),...,\gamma^{\prime}(\delta^{\prime}(t))}%
\prod \limits_{i\in I_{(t,\gamma^{\prime}(1),...,\gamma^{\prime
}(\delta^{\prime}(t)))}}\binom{m_{i}}{\sigma_{i}(t,\gamma^{\prime}(1),...,%
\gamma^{\prime}(\delta^{\prime}(t)))}
\end{align*}
such that $\ell_{(t,\gamma(1),...,\gamma(\delta(t)))}=\ell_{(t,\gamma^{%
\prime }(1),...,\gamma^{\prime}(\delta^{\prime}(t)))}^{\prime}$ and $$\prod
\limits_{1}^{\gamma(\delta(t))}d_{n,\gamma(s)}^{p^{\ell_{(t,\gamma
(1),...,\gamma(\delta(t)))}-\sum\limits_{j=1}^{s-1}(n-\gamma(j))}}=\prod%
\limits_{i\in
I_{(t,\gamma^{\prime}(1),...,\gamma^{\prime}(\delta^{\prime}(t)))}}d_{n,i}^{%
\sigma_{i}(t,\gamma^{\prime}(1),...,\gamma ^{\prime}(\delta^{\prime}(t)))}.$$
Here $\{\gamma(1),...,\gamma(\delta(t))\}$ and $\{\gamma^{\prime}(1),...,%
\gamma^{\prime}(\delta^{\prime}(t))\}$ are partitions of $\{t+1,...,n\}$ of
consecutive and non-consecutive elements respectively. For the definition of 
$I_{(t,\gamma^{\prime}(1),...,\gamma ^{\prime}(\delta^{\prime}(t)))}$ and $%
\sigma_{i}(t,\gamma^{\prime }(1),...,\gamma^{\prime}(\delta^{\prime}(t)))$,
please see the second case in the proof bellow because they strongly depend
on the particular partition.
\end{proposition}

\proof
Please see:\newline
\url{http://www.maths.warwick.ac.uk/agt/ftp/aux/agt-4-13/full.ps.gz}\hfil$\sq$\medskip

Next we consider a lemma in the ``opposite'' direction of last Proposition.

\begin{lemma}
Let $k\leq n-i$ and $i<n$, then: {\small
$d_{n,n-i}^{\alpha_{k}p^{k}+\alpha_{0}}=\left(
\Psi(d_{n,n-i}^{\alpha_{k}p^{k}+\alpha_{0}})\right) ^{\ast}+$\newline
$\binom{\alpha_{k}}{\min(\alpha_{k},\alpha_{0})}\binom{\alpha_{0}}{%
\min(\alpha_{k},\alpha_{0})}\!\left(
\Psi(d_{n,n-i-k}^{\min(\alpha_{k},\alpha_{0})p^{k}}d_{n,n-i}^{(\alpha_{k}-%
\min(\alpha
_{k},\alpha_{0}))p^{k}+(\alpha_{0}-\min(\alpha_{k},%
\alpha_{0}))}d_{n,n-i+k}^{\min(\alpha_{k},\alpha_{0})})\!\right) ^{\ast}$}
\end{lemma}

\begin{proof}
We consider all admissible sequences in $\left( i(d_{n,n-i}^{p^{k}})\right)
^{\alpha_{k}}\left( i(d_{n,n-i})\right) ^{\alpha_{0}}$.
\end{proof}

Note that $d_{n,n-i}^{\alpha_{k}p^{k}+...+\alpha_{0}}$ can be computed by
repeated use of the formulae in the last lemma for all possible choices.

\begin{remark}
\label{m-big}We must admit that if $m(n)>>0$, then there exist many
candidates for $m^{\prime}$ and the bookkeeping described above can not be
done by hand. We believe that it is harder but safer to consider all
possible choices.
\end{remark}

Next, the algorithm which calculates Adem relations using modular invariants
is demonstrated.

\begin{proposition}
\label{Adem-Rel}Let $e_{I}\in T[n]$. The following algorithm computes $\rho
(e_{I})$ in $R[n]$.

{\rm i)}\qua Let $\Re =\{m=(m_{0},...,m_{n-1})\}$ be all solutions of $|I|=\sum
\limits_{0}^{n-1}m_{i}(p^{n}-p^{i})+\sum
\limits_{1}^{l^{\prime}}(p^{n}-p^{s_{i}}-p^{k_{i}})$. Note that $s_{i}$ and $%
k_{i}$ are uniquelly defined by lemma \ref{NumBock}. Let \c{K} be the set of all admissible  sequences $K$ such that $\mid K\mid =\mid I\mid $ and $K\le I$.  Moreover, $Q_{K}\in R[n] $ and $Q_{K}=\Psi^{-1}(d^{m}M^{\epsilon})$ for $m\in \Re$.

{\rm ii)}\qua Let $h^{I^{\prime}}=\Psi_{T}^{-1}(e_{I})$ and find $b_{I,K}$ the
coefficient of $h^{I^{\prime}}$ in $\hat{\imath}(d^{m}M^{\epsilon})$ for all
elements of $\Re$.

{\rm iii)}\qua Compute the image of $d^{m}M^{\epsilon}$ in $(R[n])^{\ast}$.

{\rm iv)}\qua Use the Kronecker product to evaluate $\rho(e_{I})$ :

Start with the first non-zero $b_{I,K_{1}}$, $\ \rho(e_{I})$ \ contains $%
a_{I,K_{1}}Q_{K_{1}}$; i.e.\ $\langle d^{K_{1}^{\prime}},%
\rho(e_{I})\rangle =a_{I,K_{1}}=b_{I,K_{1}}$. Proceed to the next sequence $K_{2}$
and use $b_{I,K_{2}}$ (whether or not is zero) and the image of $%
d^{K_{2}^{\prime}}$ to compute the coefficient $a_{I,K_{2}}$ of $Q^{K_{2}}$
in $\rho(e_{I})$. Repeat last step for all remaining sequences.
\end{proposition}

We close this work by making some remarks about evaluating $\rho(e_{I})$
using matrices introduced in section 4. Since $(e_{I})^{\ast}=h^{I^{\prime}}$
is an element of $B[n]$, one has to find all sequences $m=(m_{0},%
\cdots,m_{n-1})$ such that $d^{m}$ contains $(e_{I})^{\ast}$ as a summand.
This is equivalent to find all matrices $C$ such that $(e_{I})^{\ast}=\prod%
\limits_{t=1}^{n}h_{t}^{(\mathbf{1}\cdot C)_{t-1}}$ and then group them in
different sets such that each set corresponds to an $m$. The coefficient $%
\alpha_{^{\prime}m}$ of $Q^{^{\prime}m}$ in $\rho(e_{I})$ is a function of
the order of the set corresponding to $m$. Given $h^{I^{\prime}}$, there is
a great number of choices for $C$ depending on $I^{\prime}$ as the
interested reader can easily check and this is the reason for the high
complexity of Adem relations.

\Addresses\recd

\end{document}